\documentclass[a4paper,12pt]{article}
 \usepackage[T2A]{fontenc}
\usepackage[cp1251]{inputenc}
\usepackage[english]{babel}
\usepackage{amsfonts}
\usepackage{graphics}
\usepackage{wrapfig}
\usepackage[dvips]{graphicx}

\date{}
\usepackage[tbtags]{amsmath}
\usepackage{amsfonts,amssymb}
\usepackage[usenames]{color}

\begin{document}
\bf
\begin{center} On orthogonal systems of shifts of scaling function on local fields of positive characteristic.
\end{center}
\centerline{ Gleb Sergeevich BERDNIKOV, Iuliia Sergeevna KRUSS,}
\centerline{Sergey Fedorovich LUKOMSKII}

\centerline{Departament of Mathematic Analysis, Saratov State University, Saratov, Russia.}

\footnotetext[1]{Correspondence: evrointelligent@gmail.com, KrussUS@gmail.com, LukomskiiSF@info.sgu.ru}
\rm
{\bf Abstract:} We present a new method for constructing an orthogonal step scaling function on local fields of positive characteristic, which generates multiresolution analysis.


 \section*{1. Introduction }
   Chinese mathematicians
  H.Jiang, D.Li, and N.Jin  in the
  article \cite{JLJ} introduced the notion of multiresolution analysis (MRA) on local
  fields. For the fields  $ F^{(s)}$
  of positive characteristic $p$ they proved some properties  and
  gave an algorithm for constructing wavelets for a known scaling
  function. Using these results they constructed "Haar MRA" and corresponding
  "Haar wavelets". The problem of constructing orthogonal MRA on the field $ F^{(1)}$ is studied in detail in the works \cite{YuF1,YuF3,YuF2,SL,WP,YuFWP}.

  In \cite{LJ}  a necessary
  condition and  sufficient conditions for wavelet frame on
  local fields are given. B.Behera and Q.Jahan \cite{BJ1}
  constructed the wavelet packets associated with  MRA on local
  fields of positive characteristic. In the article \cite{BJ2}   a necessary
   and  sufficient conditions for a function $\varphi\in L^2( F^{(s)})$ under which
   it is
  a scaling function for MRA are obtained.  These conditions are
   following
  \begin{equation}\label{eq0.1}
    \sum_{k\in \mathbb N_0}|\hat\varphi(\xi+u(k))|^2=1\
  \end{equation}
  for   a.e.  $\xi$ in  unit ball
    ${\cal D}$,
  \begin{equation}\label{eq0.2}
    \lim\limits_{j\to\infty}|\hat\varphi(\mathfrak p^j\xi)|=1 \ for\  a.e.\  \xi \in
    F^{(s)},
  \end{equation}
  and there exists an integral periodic function $m_0 \in L^2(\cal
  D)$ such that
  \begin{equation}\label{eq0.3}
  \hat\varphi(\xi)=m_0(\mathfrak p\xi)\hat\varphi(\mathfrak p\xi)\  for \ a.e.\  \xi
  \in F^{(s)}
  \end{equation}
    where $\{u(k) \}$ is the set of shifts, $\mathfrak p$ is a prime element.
    B.Behera and Q.Jahan \cite{BJ3} proved also  if the translates of the scaling
    functions of two multiresolution analyses are biorthogonal, then the associated
    wavelet families are also biortho\-gonal.
  So, to construct MRA on  a local field $ F^{(s)}$ we need to
  construct an integral periodic mask $m_0$ with conditions
  (\ref{eq0.1}-\ref{eq0.3}). To solve this problem in articles \cite{JLJ}, \cite{BJ3, BJ2, BJ1, LJ} was used prime element methods developed in \cite{MT}.
      In these articles    only Haar wavelets are
  obtained. In the article  \cite{SLAV}
  an another  method to construct integral periodic masks and
  corresponding  scaling step functions that generate
  non-Haar orthogonal MRA are developed.

   However, in \cite{SLAV} only simple case of mask $m_0$ being elementary is considered, i.e. $m_0(\chi)$ is constant on cosets $(F^{(s)+}_{-1})^\bot$ and $m_0(\chi)$ takes only two values 0 and 1.
   In this article, we get rid of these restrictions and specify the method of constructing the scaling function only with the  condition:  $|\hat\varphi|$ is a step function.
      We reduce this problem to the study of some dynamical system and prove that it has a fixed point.

 \section*{2. Basic concepts}

Let $p$ be a prime number, $s\in \mathbb N$, $GF(p^s)$ -- finite field. Local field $F^{(s)}$ of positive characteristic $p$ is isomorphic (Kovalski-Pontryagin theorem \cite{GGP}) to the set of formal power series
 $$
 a=\sum_{i=k}^{\infty}{\bf a}_it^i,\ k\in \mathbb{Z},\ {\bf a}_i\in GF(p^s).
 $$

 Addition and multiplication in the field $F^{(s)}$ are defined as summ and product of such series, i.e. if

 $$
 a=\sum_{i=k}^{\infty}{\bf a}_it^i,\ b=\sum_{i=k}^{\infty}{\bf b}_it^i,
 $$
 then
 $$
 a\dot+b=\sum_{i=k}^{\infty}({\bf a}_i\dot+ {\bf b}_i)t^i,\ {\bf a}_i\dot+ {\bf b}_i=({\bf a}_i + {\bf b}_i){\rm \ mod} \ p,
$$
$$
 ab=\sum_{l=2k}^{\infty}t^l\ \sum_{i,j:i+j=l}({\bf a}_i {\bf b}_j)
 $$

 Topology in $F^{(s)}$ is defined by the base of neighborhoods of zero
   $$
   F^{(s)}_n=\{a=\sum_{j=n}^\infty {\bf a}_jt^j|{\bf a}_j\in GF(p^s)\}.
   $$

 If
   $$
   a=\sum_{j=n}^\infty {\bf a}_jt^j,\ {\bf a}_n\ne {\bf 0},
   $$
  then by definition $\|a\|=(\frac{1}{p^s})^n$ which implies

   $$
   F^{(s)}_n=\{x\in F^{(s)}:\| x\|\le (\frac{1}{p^s})^n \}
   $$

  Thus we may consider local field $F^{(s)}$ of positive
  characteristic $p$ as the field of sequences infinite in both directions
  $$
  a=(\dots ,{\bf 0}_{n-1},{\bf a}_n,\dots,{\bf a}_0,{\bf a}_1,\dots),\ {\bf a}_j\in GF(p^s)
  $$
which have only finite number of elements ${\bf a}_j$ with negative $j$
  nonequal to zero, and the operations of addition and multiplication are defined by equalities
 $$
 a\dot+b=(({\bf a}_i\dot+ {\bf b}_i))_{i\in \mathbb Z},
 $$
 \begin{equation} \label{eq1.1}
 ab= (\sum_{i,j:i+j=l}({\bf a}_i {\bf b}_j))_{l\in \mathbb Z},
 \end{equation}
where $"\dot+"$ and $"\cdot"$ are respectively addition and multiplication in $GF(p^s)$.
 Thus

 $$
 \|a\|=\|(\dots,{\bf 0}_{n-1},{\bf a}_n,{\bf a}_{n+1},\dots)\|=(\frac{1}{p^s})^n, \
 \mbox{\rm если}\  {\bf a}_n\ne {\bf 0},
 $$

 $$
   F^{(s)}_n=\{a=({\bf a}_j)_{j\in \mathbb Z}: {\bf a}_j\in GF(p^s);\ {\bf a}_j=0,\ \forall j<n
   \}.
   $$

Let us consider $F^{(s)+}$ -- the additive group of the field $F^{(s)}$.
   Neighborhoods $F^{(s)}_n$ are compact subgroups of the group
   $F^{(s)+}$, we will denote them as $F^{(s)+}_n$. They have the following properties:

   1)$\dots\subset F^{(s)+}_1\subset F^{(s)+}_0\subset
   F^{(s)+}_{-1}\dots$

   2)$F^{(s)+}_n/ F^{(s)+}_{n+1}\cong GF(p^s)^+$ и $\sharp (F^{(s)+}_n/ F^{(s)+}_{n+1})=p^s$.

This implies that
   if $s=1$ then $F^{(1)+}$ is Vilenkin group with the stationary generating sequence $p_n=p$. The inverse is also true:  one can define multiplication in any Vilenkin group $(\mathfrak G,\dot +)$ with stationary generating sequence $p_n=p$ using equality (\ref{eq1.1}). Supplied with such operation
   $(\mathfrak G,\dot+,\cdot)$ becomes a field isomorphic to $F^{(1)}$, where
   $e=(\dots,0,0_{-1},1_0,0_1,\dots)$ is a neutral element with respect to multiplication.

It is noted in \cite{AVSL} that the field $F^{(s)}$ can be described as a linear space over $GF(p^s)$. Using this description one may define the multiplication of element $a\in F^{(s)} $ on element $\overline\lambda \in GF(p^s)$ coordinatewise, i.e. $\overline\lambda a =(\dots {\bf 0}_{n-1},\overline\lambda {\bf a}_n,\overline\lambda {\bf a}_{n+1},\dots)$, and the modulus
    $\overline\lambda \in GF(p^s)$ can be defined as
$$ |\overline\lambda|=\left\{
           \begin{array}{ll}
 1,& \overline\lambda \ne {\bf 0},\\
 0,& \overline\lambda = {\bf 0}.\\
  \end{array} \right.
$$

 It is also proved there, that the system $g_k\in F_k^{(s)}\setminus F_{k+1}^{(s)}$ is a basis in $F^{(s)}$, i.e. any element $a\in F^{(s)}$ can be represented as:
    \begin{center}
      $a=\sum\limits_{k\in\mathbb{Z}}\overline\lambda_kg_k,\
      \overline\lambda_k\in GF(p^s)$. \\
    \end{center}

    From now on we will consider $g_k=(...,{\bf 0}_{k-1},(1^{(0)},0^{(1)},...,0^{(s-1)})_k,{\bf 0}_{k+1},...)$. In this case $\overline\lambda_k={\bf a}_k$.

 Let us define the sets
 $$
    H_0^{(s)}=\{h\in G: h={\bf a}_{-1}g_{-1}\dot+{\bf a}_{-2}g_{-2}\dot+\dots \dot+ {\bf a}_{-s}g_{-s}
    \},s\in \mathbb N.
  $$
$$
H_0=\{h\in
G:\;h={\bf a}_{-1}g_{-1}\dot+{\bf a}_{-2}g_{-2}\dot+\dots\dot+{\bf a}_{-s}g_{-s},\;s\in\mathbb
N\}.
$$

 The set $H_0$ is the set of shifts in $F^{(s)}$. It is an analogue of the set of nonnegative integers.

 We will denote the collection of all characters of $F^{(s)+}$ as $X$. The set $X$ generates a commutative group with respect to the multiplication of characters: $(\chi*\phi) (a)=\chi(a)\cdot \phi(a)$. Inverse element is defined as $\chi^{-1}(a)=\overline{\chi(a)}$, and the neutral element is $e(a)\equiv1$.

 Following [15] we define characters $r_n$ of the group $F^{(s)+}$ in the following way.
 Let $x=(\dots,{\bf 0}_{k-1},{\bf x}_k,$ ${\bf x}_{k+1},\dots)$, ${\bf x}_j=(x_j^{(0)},x_j^{(1)},\dots,x_j^{(s-1)})\in GF(p^s)$. The element ${\bf x}_j$ can be written in the form ${\bf x}_j=(x_{js+0},x_{js+1},\dots,x_{js+(s-1)})$. In this case
 $$ x=(\dots,0,...,0,x_{ks+0},x_{ks+1},\dots,x_{ks+s-1},x_{(k+1)s+0},x_{(k+1)s+1},\dots,x_{(k+1)s+s-1},\dots)
 $$
 and the collection of all such sequences $x$ is Vilenkin group. Thus the equality
 $r_n(x)=r_{ks+l}(x)=e^{\frac{2\pi i}{p}(x_{ks+l})}$ defines Rademacher function of $F^{(s)+}$ and every character $\chi\in X$ can be described in the following way:

\begin{equation} \label{eq1.2}
 \chi=\prod \limits_{n\in\mathbb Z}{r}_n^{ a_n},\quad a_n=\overline{0,p-1}.
\end{equation}
The equality (\ref{eq1.2}) can be rewritten as
\begin{equation} \label{eq1.3}
  \chi=
  \prod\limits_{k\in\mathbb Z} r_{ks+0}^{a_k^{(0)}}r_{ks+1}^{a_k^{(1)}}\dots r_{ks+s-1}^{a_k^{(s-1)}}
 \end{equation}
 and let us define
   $$
   r_{ks+0}^{a_k^{(0)}}r_{ks+1}^{a_k^{(1)}}\dots r_{ks+s-1}^{a_k^{(s-1)}}={\bf r}_k^{{\bf a}_k}
   $$
   where ${\bf a}_k=(a_k^{(0)},a_k^{(1)},\dots,a_k^{(s-1)})\in GF(p^s)$. Then (\ref{eq1.3}) takes the form
    \begin{equation} \label{eq1.4}
     \chi =\prod_{k\in \mathbb Z}{\bf r}_k^{{\bf a}_k}.
    \end{equation}

    We will refer to ${\bf r}_k^{(1,0,\dots,0)}={\bf r}_k$ as the Rademacher functions.
    By definition we set

    $$
   ({\bf r}_k^{{\bf a}_k})^{{\bf b}_k}={\bf r}_k^{{\bf a}_k{{\bf b}_k}}, \quad
     \chi^{\bf b}=(\prod {\bf r}_k^{{\bf a}_k})^{\bf b}=\prod {\bf r}_k^{{\bf a}_k\bf b},
     \quad {\bf a}_k, {\bf b}_k, {\bf b}\in
  GF(p^s).
          $$

The definition of Rademacher function implies that if ${\bf x}=((x_k^{(0)},x_k^{(1)},\dots x_k^{(s-1)}))_{k\in \mathbb Z}$ and ${\bf u}=(u^{(0)},u^{(1)},\dots, u^{(s-1)})\in GF(p^s)$ then
  $$
({\bf r}_k^{{\bf u}},{\bf x})=\prod\limits_{l=0}^{s-1}e^{\frac{2\pi
i}{p}u^{(l)}x_k^{(l)}}.
$$

In [15] the following properties of characters are proved

1) ${\bf r}_k^{{{\bf u}\dot+{\bf v}}}={\bf r}_k^{\bf u}{\bf r}_k^{\bf v}$, ${\bf u}, {\bf v}\in GF(p^s)$.

2) $({\bf r}_k^{\bf v},{\bf u}g_j)=1$, $\forall k\ne j$, ${\bf u}, {\bf v}\in GF(p^s)$.

3)  The set of characters of the field $F^{(s)}$ is a linear space
   $(X,\; *,\; \cdot^{GF(p^s)})$
   over the finite field $GF(p^s)$ with multiplication being an inner operation and the power ${\bf u}\in GF(p^s)$being an outer operation.

4) The sequence of Rademacher functions $({\bf r}_k)$
   is a basis in the space $(X,\; *,\; \cdot^{GF(p^s)})$.

5) Any sequence of characters $\chi_k\in (F_{k+1}^{(s)})^\bot\setminus(F_k^{(s)})^\bot$
   is also a basis in the space $(X,\; *,\; \cdot^{GF(p^s)})$, where ${F^{(s)}_n}^{\bot}$ is the annihilator of $F^{(s)+}_n$.

The dilation operator
  ${\cal A}$ in local field $F^{(s)}$ can be defined as ${\cal
A}x:=\sum_{n=-\infty}^{+\infty}{\bf a}_ng_{n-1}$, where
 $x=\sum_{n=-\infty}^{+\infty}{\bf a}_ng_n\in F^{(s)}$. In the group of characters it is defined as
 $(\chi {\cal A},x)=(\chi, {\cal
 A}x)$.

 \section*{3. Scaling function and MRA}

  We will consider a case of scaling function $\varphi$, which generates an orthogonal MRA, being step function. The set of step functions constant on cosets of a subgroup $F_M^{(s)}$
  with the support ${\rm supp}(\varphi)\subset F^{(s)}_{-N}$ will be denoted as
  $\mathfrak D_M(F^{(s)}_{-N})$, $M,N\in \mathbb N$. Similarly, $\mathfrak D_{-N}({F^{(s)}_{M}}^\bot)$ is a set of step functions, constant on the cosets of a subgroup ${F^{(s)}_{-N}}^{\bot}$
  with the support ${\rm supp}(\varphi)\subset {F^{(s)}_M}^\bot$. If $\varphi\in \mathfrak D_M(F^{(s)}_{-N})$
  generates an orthogonal MRA, it satisfies the refinement equation
   $\varphi(x)=\sum_{h\in H_0^{(N+1)}}\beta_h\varphi({\cal
   A}x\dot-h)$\ \cite{SLAV}, which can be rewritten in a frequency from

    \begin{equation}                                      \label{eq2.1}
 \hat\varphi(\chi)=m_0(\chi)\hat\varphi(\chi{\cal
  A}^{-1}),
 \end{equation}
 where
 \begin{equation}                                           \label{eq2.2}
 m_0(\chi)=\frac{1}{p}\sum_{h\in
 H_0^{(N+1)}}\beta_h\overline{(\chi{\cal A}^{-1},h)}
 \end{equation}

 is the mask of equation (\ref{eq2.1}).

For the step functions in \cite{SLAV}  condition (\ref{eq0.3}) and orthogonality condition (\ref{eq0.1}) are rewritten in the terms of Rademacher functions

1) If $\hat\varphi(\chi)\in\mathfrak
 D_{-N}({F^{(s)}_M}^\bot)$ is
  a solution of refinement equation (\ref{eq2.1}) and the system
  of shifts $(\varphi(x\dot-h))_{h\in H_0}$ is orthonormal,
  then $\varphi$ generates an orthogonal MRA.

 2) If
  $\hat\varphi(\chi)\in\mathfrak
 D_{-N}({F^{(s)}_M}^\bot)$
 , then the system of shifts
 $(\varphi(x\dot-h))_{h\in H_0}$ will be orthonormal
  iff for any
${\bf a}_{-N},{\bf a}_{-N+1},\dots,{\bf a}_{-1}\in GF(p^s)$
 \begin{equation}                                                \label{eq2.3}
 \sum_{{\bf a}_{0},{\bf a}_1,\dots,{\bf a}_{M-1}\in GF(p^s)}|\hat\varphi({F^{(s)}_{-N}}^\bot
 {\bf r}_{-N}^{{\bf a}_{-N}}\dots {\bf r}_0^{{\bf a}_0}\dots
 {\bf r}_{M-1}^{{\bf a}_{M-1}})|^2=1.
 \end{equation}

Thus to construct an orthogonal MRA one must construct a function $\hat\varphi(\chi)\in\mathfrak
   D_{-N}({F^{(s)}_M}^\bot)$, which is
  a solution of refinement equation (\ref{eq2.1}) and which satisfies conditions (\ref{eq2.3}). Satisfying both conditions is the main difficulty of this problem.

  As it was already mentioned in introduction, a method for construction of scaling function which generates nonhaar orthogonal MRA is specified in \cite{SLAV}. It is constructed by the means of some tree and results in a function such that $|\varphi|$ takes two values only: 0 and 1. More general case will be presented in the next section.

  \section*{4. Construction of orthogonal scaling function}

{\bf Definition 4.1.} {\it  Let $F^{(s)}$ be a local field of positive characteristic $p$, $N$ is a natural number. Then by $N$-valid tree we mean a tree,
oriented from leaves to root and satisfying conditions:

 1)Every vertex is an element of $GF(p^s)$, i.e has the form  ${\bf a}_{i}=(a_i^{(0)},a_i^{(1)},\dots,a_i^{(s-1)})$, $a_i^{(j)}=\overline{0,p-1}$.

 2)The root and all vertices of level $N-1$ are equal to the zero element of $GF(p^s)$: ${\bf 0}=(0^{(0)},0^{(1)},\dots,0^{(s-1)})$.

 3)Any path $({\bf a}_k\to{\bf a}_{k+1}\to\dots\to{\bf a}_{k+N-1})$ of length $N-1$
  appears in the tree exactly one time.}

  Let us choose  $N$-valid tree $T$ and construct a scaling function using it.

1) We will use this tree $T$ to construct new tree $\tilde T$. Every vertex of the tree $\tilde T$ is a vector of $N$ elements each being an element of $GF(p^s)$: ${\bf A}=({\bf a}_N,{\bf a}_{N-1},\dots,{\bf a}_1)$. Such vertices are constructed in the following way: if a tree $T$ has a path of length $N-1$ starting from ${\bf a}_N$
 $$
 {\bf a}_N\rightarrow{\bf a}_{N-1}\rightarrow\dots\rightarrow{\bf a}_1,
 $$
 then in $\tilde T$ we will have a vertex with the value
 equal to the array of $N$ elements $({\bf a}_N,{\bf a}_{N-1},\dots$ $\dots,{\bf a}_1)$.  Due to condition 3) of $N$-validity of tree $T$ each such array corresponds to the unique vertex of the new tree $\tilde T$. Thus, the root of $\tilde T$ is an $N$-dimensional vector with all elements equal to zero of $GF(p^s)$ ${\bf O}=({\bf 0},{\bf 0},\dots,{\bf 0})$. Vertices of level 1 in the tree $\tilde T$ are $N$-dimensional vectors, which have all their elements, except the first one, equal to zero of $GF(p^s)$: $({\bf a}_i,{\bf 0},\dots,{\bf 0})$, where ${\bf a}_i$ is some vertex of level $N$ in the tree $T$. Vertices of level 2 in the tree $\tilde T$ are $N$-dimensional vectors: $({\bf a}_{i_2},{\bf a}_{i_1},{\bf 0},\dots,{\bf 0})$, where ${\bf a}_{i_2}$ and ${\bf a}_{i_1}$ are some vertices of levels $N+1$ and $N$ of the tree $T$ respectively, which are connected. We should note that in this example ${\bf a}_{i_1}\ne{\bf 0}$, but ${\bf a}_{i_2}$ may be zero element of $GF(p^s)$. Thus in $\tilde T$ connected vertices have the form: $({\bf a}_{i_N},{\bf a}_{i_{N-1}},\dots,{\bf a}_{i_1})\rightarrow({\bf a}_{i_{N-1}},
 \dots,{\bf a}_{i_1},{\bf a}_{i_0})$. However not all vertices satisfying this condition will be connected. Arcs are taken from the original tree $T$.
 If we denote $height(T)=H,$ $height(\tilde T)=\tilde H$,
 then obviously $\tilde H=H-N+1$.

2) Now we will construct a directed graph $\Gamma$ using $\tilde T$.
 We connect each vertex
 ${\bf A}_N=({\bf a}_N,{\bf a}_{N-1},\dots,{\bf a}_1)$ of
 $\tilde T$  to each vertex of lesser level of the form
 $({\bf a}_{N-1},\dots,{\bf a}_1,{\bf a}_0)$, i.e having first $(N-1)$ elements equal to the last $(N-1)$ elements of vertex ${\bf A}_N$.
 The vertices, to which
 ${\bf A}_N$ is connected, we will denote by
 $({\bf a}_{N-1},\dots,{\bf a}_1,\tilde {\bf a}_0)$. I.e.
 ${\bf a}_0\in \{\tilde {\bf a}_0\}$ iff the vertex ${\bf A}_N$ is connected to
 $({\bf a}_{N-1},\dots,{\bf a}_1,{\bf a}_0)$ in digraph $\Gamma$.

3)
Let us denote

 $$
 \lambda_{{\bf a}_{-N},{\bf a}_{-N+1},\dots,{\bf a}_{-1},{\bf a}_0}=
 |m_0({F^{(s)}}^\bot_{-N}{\bf r}_{-N}^{{\bf a}_{-N}}{\bf r}_{-N+1}^{{\bf a}_{-N+1}}\dots
 {\bf r}_{-1}^{{\bf a}_{-1}}{\bf r}_{0}^{{\bf a}_{0}})|^2,
 $$
 i.e.
 $\lambda_{{\bf a}_{-N},{\bf a}_{-N+1},\dots,{\bf a}_{-1},{\bf a}_0} $
 is an $(N+1)$-dimensional array, enumerated by the elements of $GF(p^s).$

  If the vertex
$({\bf a}_{N},{\bf a}_{N-1},\dots,{\bf a}_{1})$ of graph $\Gamma$
 is connected to the vertices $({\bf a}_{N-1},{\bf a}_{N-2}\dots,{\bf a}_{1},$ $\tilde{\bf a}_0)$
 then we define the values of the mask in the way satisfying the condition
 \begin{equation}\label{eq3.1}
 \sum\limits_{\tilde{\bf a}_0}
 \lambda_{{\bf a}_{-N},{\bf a}_{-N+1},\dots,{\bf a}_{-1},\tilde{\bf a}_0}=1
  \ \mbox{and}\
 \lambda_{{\bf a}_{-N},{\bf a}_{-N+1},\dots,{\bf a}_{-1},{\bf a}_0}=0\  \mbox{for any ${\bf a}_0\notin\{\tilde{\bf a}_0\}.$}
 \end{equation}
  Also, let us define
 $m_0({F^{(s)}_{-N}}^\bot)=1,$ which implies
 $\lambda_{{\bf 0},{\bf 0},\dots,{\bf 0}}=1$.

To present the main result we will need some extra notation.
 Firsly, we must note that the orthonormality condition (\ref{eq2.3}) for the system of shifts of
 $\varphi(x)$ can be rewritten as: for any ${\bf a}_{-N},{\bf a}_{-N+1},\dots,{\bf a}_{-1}\in GF(p^s)$
  $$
  1=\sum_{{\bf a}_0,{\bf a}_1,\dots,{\bf a}_{M-1}\in GF(p^s)}|\hat\varphi(
  {F^{(s)}_{-N}}^\bot {\bf r}_{-N}^{{\bf a}_{-N}}\dots
  {\bf r}_{-1}^{{\bf a}_{-1}}{\bf r}_{0}^{{\bf a}_0}\dots
  {\bf r}_{M-1}^{{\bf a}_{M-1}})|^2=
   $$
  $$ =\sum_{{\bf a}_0\in GF(p^s)}\lambda_{{\bf a}_{-N},{\bf a}_{-N+1},\dots,{\bf a}_0}
  \sum_{{\bf a}_1\in GF(p^s)}\lambda_{{\bf a}_{-N+1},{\bf a}_{-N+2},\dots,{\bf a}_1}\dots
  $$
  $$
  \dots\sum_{{\bf a}_{M-2}\in GF(p^s)}\lambda_{{\bf a}_{M-N-2},{\bf a}_{M-N-1},\dots,{\bf a}_{M-2}}
   $$

 \begin{equation}                                                    \label{eq3.2}
  \sum_{{\bf a}_{M-1}\in GF(p^s)}\lambda_{{\bf a}_{M-N-1},{\bf a}_{M-N},\dots,{\bf a}_{M-1}}\lambda_{{\bf a}_{M-N},
  {\bf a}_{M-N+1},
  \dots,{\bf a}_{M-1},{\bf 0}}
   \dots\lambda_{{\bf a}_{M-1},{\bf 0},\dots,{\bf 0}}.
  \end{equation}

  Let us then define a sequence of  $N$-dimensional arrays
  $A^{(n)}=(a^{(n)}_{{\bf i}_1, {\bf i}_2,\dots, {\bf i}_N})_{{\bf i}_1,{\bf i}_2,\dots,{\bf i}_N\in GF(p^s)}$
  recurrently by giving the relations of their components:

  \begin{equation}
  \label{eq3.3}
  a^{(0)}_{{\bf i}_1,{\bf i}_2,\dots,{\bf i}_N}=\lambda_{{\bf i}_1,{\bf i}_2,\dots,{\bf i}_N,{\bf 0}}\lambda_{{\bf i}_2,{\bf i}_3,\dots,{\bf i}_N,{\bf 0},{\bf 0}}
  \dots\lambda_{{\bf i}_N,{\bf 0},\dots,{\bf 0}},
  \end{equation}
  \begin{equation}
  \label{eq3.4}
  a^{(n)}_{{\bf i}_1,{\bf i}_2,\dots,{\bf i}_N}=\sum_{{\bf j}\in GF(p^s)}\lambda_{{\bf i}_1,{\bf i}_2,\dots,{\bf i}_N,{\bf j}}a^{(n-1)}_{{\bf i}_2,{\bf i}_3,\dots,{\bf i}_N,{\bf j}}
  \end{equation}
We will say that the element
$a^{(s)}_{{\bf i}_1,{\bf i}_2,\dots,{\bf i}_N}$ corresponds to vertex
$({\bf i}_1,{\bf i}_2,\dots,{\bf i}_N)$.

Using new notation, orthonormality condition (\ref{eq3.2}) can be reformulated in the following way:
 the system of shifts of the function $\varphi(x)\in \mathfrak
D_{M}(F^{(s)}_{-N})$ is orthonormal if and only if
 for any ${\bf i}_1,{\bf i}_2,\dots,{\bf i}_N$: $a^{(M)}_{{\bf i}_1,{\bf i}_2,\dots,{\bf i}_N}=1$, in other words, iff an array $A^{(M)}$ has all its elements equal to 1.

{\bf Lemma 4.1.} {\it  The components of $A^{(0)}$ corresponding to vertices of level
$l\leq N$ in the tree $\tilde T$ are equal to 1.}

{\bf Proof.}
Firstly, let us notice that any vertex of $\tilde T$ of level
$l\leq N$ has the form
$({\bf a}_l,{\bf a}_{l-1},\dots,{\bf a}_1,{\bf 0},\dots,{\bf 0}),\quad
{\bf a}_1\neq {\bf 0}$.  Indeed,
if a vertex has level $l$ in $\tilde T$, then the first element of the vector - the vertex of $T$ - is of level $l+N-1$ in $T$ and is the beginning of the following path directed to root:
$({\bf a}_l\rightarrow{\bf a}_{l-1}\rightarrow\dots\rightarrow{\bf a}_1\rightarrow
{\bf 0}\rightarrow\dots\rightarrow {\bf 0})$, where ${\bf a}_1$ is a vertex of level
$N$ and is nonzero by the $N$-validity condition.

We will prove {\it the lemma} by induction on $l$. Let
$l=0$. Thus, we consider the root of $\tilde T$. The root has the form
$({\bf 0},{\bf 0},\dots,{\bf 0})$. By construction $\lambda_{{\bf 0},{\bf 0},\dots,{\bf 0}}=1$. Its corresponding element of array $A^{(0)}$ is $a^{(0)}_{{\bf 0},{\bf 0},\dots,{\bf 0}}$.
Let us substitute ${\bf i}_1,{\bf i}_2,\dots,{\bf i}_N={\bf 0}$ into (\ref{eq3.3}). We obtain
$$a^{(0)}_{{\bf 0},{\bf 0},\dots,{\bf 0}}=\lambda_{{\bf 0},{\bf 0},\dots,{\bf 0}}\lambda_{{\bf 0},{\bf 0},\dots,{\bf 0}}\dots\lambda_{{\bf 0},{\bf 0},\dots,{\bf 0}}=1.$$

Now we prove if any vertex of level $l= k-1< N$
satisfies the condition
$a^{(0)}_{{\bf a}_{k-1},{\bf a}_{k-2},\dots,{\bf a}_1,{\bf 0},\dots,{\bf 0}}=1$, then
such condition is also satisfied by any vertex of level $l=k\leq N$ of the tree
$\tilde T$. Using (\ref{eq3.3}) and substituting
${\bf i}_1={\bf a}_{k-1},{\bf i}_2={\bf a}_{k-2},\dots,{\bf i}_{k-1}={\bf a}_1\neq
{\bf 0},{\bf i}_{k}={\bf 0},\dots,{\bf i}_N={\bf 0},$ we rewrite the induction hypothesis:

$$
a^{(0)}_{{\bf a}_{k-1},{\bf a}_{k-2},\dots,{\bf a}_1,{\bf 0},\dots,{\bf 0}}=\lambda_{{\bf a}_{k-1},{\bf a}_{k-2},\dots,{\bf a}_1,{\bf 0},\dots,{\bf 0}}
\lambda_{{\bf a}_{k-2},{\bf a}_{k-3},\dots,{\bf a}_1,{\bf 0},\dots,{\bf 0}}
\dots\lambda_{{\bf a}_1,{\bf 0},\dots,{\bf 0}}\lambda_{{\bf 0},{\bf 0},\dots,{\bf 0}}\dots
$$
\begin{align*}
\label{eq4}
\dots\lambda_{{\bf 0},{\bf 0},\dots,{\bf 0}}=\lambda_{{\bf a}_{k-1},{\bf a}_{k-2},\dots,{\bf a}_1,{\bf 0},\dots,{\bf 0}}\lambda_{{\bf a}_{k-2},{\bf a}_{k-3},\dots,{\bf a}_1,{\bf 0},\dots,{\bf 0}}
\dots\lambda_{{\bf a}_1,{\bf 0},\dots,{\bf 0}}=1
\end{align*}
Here we omit $\lambda_{{\bf 0},{\bf 0},\dots,{\bf 0}}=1$. Now, let
$$
{\bf A}_k=({\bf a}_k,{\bf a}_{k-1},\dots,{\bf a}_1,{\bf 0},\dots,{\bf 0}),\quad
{\bf a}_1\neq {\bf 0}
$$
 be a vertex of level $k$ of $\tilde T$.

Let this vertex be connected to the vertex ${\bf A}_{k-1}=({\bf a}_{k-1},\dots,{\bf a}_1,{\bf 0},\dots,{\bf 0})$
of level $k-1$ in $\tilde T$. Then it can be shown that the vertex ${\bf A}_k$ is only connected to the vertex
${\bf A}_{k-1}$ in digraph $\Gamma$ also.

Firstly, let us prove that in graph $\Gamma$ the vertex
${\bf A}_k$ is not connected to any other vertex, which has level $k-1$ in $\tilde T$. We will prove the fact by contradiction.
 Assume that
${\bf B}_{k-1}=({\bf b}_{k-1},\dots,{\bf b}_1\neq{\bf 0},{\bf 0},\dots,{\bf 0})$
is another vertex which has level $k-1$ in $\tilde T$ and that
${\bf A}_k$ is connected to ${\bf A}_{k-1}$ and
${\bf B}_{k-1}$ in graph $\Gamma$. By construction, if
${\bf A}_k$ is connected to ${\bf B}_{k-1}$ then for any
$i=\overline{1,k-1},\quad {\bf a}_i={\bf b}_i$, which implies
vertices ${\bf A}_{k-1}$ and ${\bf B}_{k-1}$ being identical,
which contradicts the uniqueness of the vertices in $\tilde
T$ and $\Gamma$. Thus, there is only one vertex, which is of level
$(k-1)$ in $\tilde T$ and to which
${\bf A}_k$ is connected in graph $\Gamma$.

Secondly, we prove that
in $\Gamma$ the vertex ${\bf A}_k$ is not connected to any vertex, which has level strictly less, than $k-1$ in the tree $\tilde T$.
 Let $n>1$,
${\bf B}_{k-n}=({\bf b}_{k-n},\dots,{\bf b}_1,{\bf 0},\dots,{\bf 0})$ be
an arbitrary vertex of level $(k-n)$ in $\tilde T$. By construction of $\Gamma$,
for the vertex ${\bf A}_k$ to be connected to ${\bf B}_{k-n}$ it is necessary for the equality
${\bf a}_1={\bf 0}$ to hold, which is impossible by assumption
${\bf a}_1\neq{\bf 0}$. Thus, we proved that the vertex
${\bf A}_k$ is connected only to
${\bf A}_{k-1}$ in $\Gamma$.

By construction that means that $\lambda_{{\bf a}_k,\dots,{\bf a}_1,{\bf 0},\dots,{\bf 0}}=1$.
Thus, substituting
${\bf i}_1={\bf a}_k,{\bf i}_2={\bf a}_{k-1},\dots,{\bf i}_k={\bf a}_1,{\bf i}_{k+1}={\bf 0},\dots,{\bf i}_N={\bf 0}$
into (\ref{eq3.3}) and using the induction hypothesis
we obtain

$$a^{(0)}_{{\bf a}_k,\dots,{\bf a}_1,{\bf 0},\dots,{\bf 0}}=\lambda_{{\bf a}_k,\dots,{\bf a}_1,{\bf 0},\dots,{\bf 0}}
\lambda_{{\bf a}_{k-1},\dots,{\bf a}_1,{\bf 0},\dots,{\bf 0}}\dots\lambda_{{\bf a}_1,{\bf 0},\dots,{\bf 0}}=$$
$$=\lambda_{{\bf a}_k,\dots,{\bf a}_1,{\bf 0},\dots,{\bf 0}}a^{(0)}_{{\bf a}_{k-1},\dots,{\bf a}_1,{\bf 0},\dots,{\bf 0}}=1.$$

Lemma is proved.

{\bf Lemma 4.2.} {\it Let us consider $N$-valid tree $T$ and tree $\tilde T$ and
digraph $\Gamma$ constructed using it. Let the values of $m_0(\chi)$ be defined as specified in equalities (\ref{eq2.1}). Let also $(A^{(n)})_{n=0}^\infty$ be
a sequence of arrays defined by equalities
(\ref{eq3.3}) and ($\ref{eq3.4}$). Then the array $A^{(n)}$
has its elements corresponding to the vertices of level $l\leq N+n$ in the tree
$\tilde T$ equal to 1.}

{\bf Proof.} We will prove the lemma by induction.
The validity of base for $n=0$ follows from the previous lemma.
Now we prove that if in $A^{(n-1)}$ elements corresponding to vertices of level less or equal to $N+n-1$ are equal to one, then in $A^{(n)}$
elements corresponding to vertices of level less or equal to $N+n$ are equal to one. Let
${\bf A}_N=({\bf a}_{N},{\bf a}_{N-1},\dots,{\bf a}_1)$ be
a vertex of level $l\leq N+n$ in $\tilde T$. In graph $\Gamma$
it is connected to al vertices of lower level, which we denote as $({\bf a}_{N-1},\dots,{\bf a}_1,\tilde{\bf a}_0)$, moreover
$\sum\limits_{\tilde{\bf a}_0}
\lambda_{{\bf a}_{N},{\bf a}_{N-1},\dots,{\bf a}_{1},\tilde{\bf a}_0}=1$
and $\lambda_{{\bf a}_{N},{\bf a}_{N-1},\dots,{\bf a}_{1},{\bf a}_0}=0 \
\forall {\bf a}_0\notin\{\tilde{\bf a}_0\}.$

Also, it should be mentioned that since vertices
$({\bf a}_{N-1},\dots,{\bf a}_1,\tilde{\bf a}_0)$ of $\tilde T$
have their level not higher than $l-1\leq N+n-1$, then, by the induction hypothesis

$a^{(n-1)}_{{\bf a}_{N-1},\dots,{\bf a}_1,\tilde{\bf a}_0}=1$,
$\forall \tilde{\bf a}_0\in\{\tilde{\bf a}_0\}.$
 Then
$$a^{(n)}_{{\bf a}_{N},{\bf a}_{N-1}\dots,{\bf a}_1}=\sum_{{\bf a}_0\in GF(p^s)} \lambda_{{\bf a}_{N},{\bf a}_{N-1},\dots,{\bf a}_1,{\bf a}_0}a^{(n-1)}_{{\bf a}_{N-1},\dots,{\bf a}_1,{\bf a}_0}=$$
$$=\sum_{\tilde{\bf a}_0\in\{\tilde{\bf a}_0\}}\lambda_{{\bf a}_{N},{\bf a}_{N-1},\dots,{\bf a}_1,\tilde{\bf a}_0}
a^{(n-1)}_{{\bf a}_{N-1},\dots,{\bf a}_1,\tilde{\bf a}_0}=\sum_{\tilde{\bf a}_0\in\{\tilde{\bf a}_0\}}\lambda_{{\bf a}_{N},{\bf a}_{N-1},\dots,{\bf a}_1,\tilde{\bf a}_0}=1$$
which proves the lemma.\\
 These lemmas directly imply the following theorem.

{\bf Theorem 4.3.} {\it Let the tree $\tilde T$ and digraph $\Gamma$ be constructed using
 $N$-valid tree $T$. Let the values of $m_0(\chi)$ be defined as specified by equalities (\ref{eq3.1}). Let $\tilde H=height(\tilde T).$
 Then the equality
 $$
 \hat\varphi(\chi)=\prod\limits_{k=0}^\infty m_0(\chi{\cal
 A}^{-k})\in \mathfrak D_{-N}({F^{(s)}_M}^\bot)
 $$
 defines an orthogonal scaling function
 $\varphi(x)\in\mathfrak D_M({F^{(s)}_{-N}})$, and $M\leq\tilde
 H-N$.}

{\bf Remark.} Let us denote
the collection of functions $f_N:\{0,1,...,p-1\}^N\to [0,1]$ as $\Phi_N$
and choose a function $\Lambda\in \Phi_{N+1}$.  Function $\Lambda$ may be viewed as $N+1$-dimensional array $\Lambda=
(\lambda_{i_1,i_2,...,i_N,i_{N+1}})$. Then the equalities
(\ref{eq3.4}) define discrete dynamic system
$\Lambda:\Phi_N\to\Phi_N$, and the equality (\ref{eq3.3}) defines
the initial trajectory point. Theorem 4.3 specifies a class of discrete systems $\Lambda$, which have a stationary point in their trajectory starting from initial point (\ref{eq3.3}).

The theorem 4.3 for $s=1,\ N=1$ was proved by Iu.Kruss, for $s=1,\ N\in \mathbb N $ -- by G.Berdnikov, for any $s, N\in \mathbb N$ -- by Iu.Kruss. The idea to consider local field of positive characteristic as vector space was proposed by S.Lukomskii.

The results were obtained within the framework of the state task of
Russian Ministry of Education and Science (project 1.1520.2014K).

\end{document}